\documentclass[secthm,seceqn,number,dvips]{arxbj}
\usepackage{mathbh}

\aid{0}
\volume{13}
\issue{2}
\pubyear{2007}
\firstpage{473}
\lastpage{491}
\doi{10.3150/07-BEJ5123}

\makeatletter

\newcommand{\Xgen}{X}

\newcommand{\Nset}{\mathbb{N}}

\newcommand{\Rset}{\mathbb{R}}
\newcommand{\Zset}{\mathbb{Z}}

\newcommand{\esp}{\mathbb{E}}
\newcommand{\pr}{\mathbb{P}}
\newcommand{\XS}{X_S}
\newcommand{\mcl}{\mathcal{L}}

\newcommand{\var}{\operatorname{var}}
\newcommand{\cov}{\operatorname{cov}}
\newcommand{\cum}{\mathrm{cum}}
\newcommand{\Leb}{\mathrm{Leb}}
\newcommand{\bd}{\mathbf{d}}
\renewcommand{\d}{\,\mathrm{d}}
\newcommand{\ind}{\mathbh{1}}
\newcommand{\mcn}{\mathcal{N}}
\newcommand{\mce}{\mathcal{E}}
\newcommand{\NS}{N_S}

\newcommand{\interp}{\mathrm{I}_{\phi}}
\newcommand{\interpbis}{\overline{\mathrm{I}}_{\phi}}
\newcommand{\nj}{n}
\newcommand{\rme}{\mathrm{e}}
\newcommand{\rmi}{\mathrm{i}}
\newcommand{\mcs}{\mathcal{S}}
\newcommand{\eqref}[1]{(\ref{#1})}

\newtheorem{hypo}{Assumption}
\newtheorem{prop}[thm]{Proposition}
\newtheorem{lem}[thm]{Lemma}

\newremark{rem}{Remark}[section]
\newremark{xmpl}{Example}[section]

\makeatother

\begin{document}
\begin{frontmatter}

\title{Estimation of the memory parameter of the infinite-source Poisson process}
\runtitle{The infinite-source Poisson process}

\begin{aug}

\author[a]{\fnms{Gilles} \snm{Fa\"y}\thanksref{a}\ead[label=e1]{gilles.fay@univ-lille1.fr}\corref{}},
\author[b]{\fnms{Fran\c{c}ois} \snm{Roueff}\thanksref{b}\ead[label=e3]{roueff@tsi.enst.fr}}
\and
\author[c]{\fnms{Philippe} \snm{Soulier}\thanksref{c}\ead[label=e2]{philippe.soulier@u-paris10.fr}}

\pdfauthor{Gilles Fay, Francois Roueff, Philippe Soulier}

\runauthor{G. Fa\"y, F. Roueff and P. Soulier}

\address[a]{Laboratoire Paul-Painlev\'e, Universit\'e de Lille 1, 59655 Villeneuve
  d'Ascq Cedex, France.\\ \printead{e1}}

\address[b]{\'Ecole Nationale Sup\'erieure des T\'el\'ecommunications, 46, rue Barrault,
75634 Paris Cedex 13, France. \printead{e3}}

\address[c]{Modal'X, Universit\'e de Paris X, 200 avenue de la R\'epublique,
92000 Nanterre, France.\\ \printead{e2}}

\end{aug}

\received{\smonth{9} \syear{2005}}
\revised{\smonth{11} \syear{2006}}

\begin{abstract}
Long-range dependence induced by heavy tails is a widely reported feature of
internet traffic. Long-range dependence can be defined as the regular
variation of the variance of the integrated process, and half the index of
regular variation is then referred to as the Hurst index.
The infinite-source Poisson process (a particular case of which is the $M/G/\infty$ queue)
is a simple and popular model with this property, when the tail of the
service time distribution is regularly varying. The Hurst index of the
infinite-source Poisson process is then related to the index of regular
variation of the service times. In this paper, we present a wavelet-based
estimator of the Hurst index of this process, when it is observed either
continuously or discretely over an increasing time interval. Our estimator is
shown to be consistent and robust to some form of non-stationarity. Its rate
of convergence is investigated.
\end{abstract}


\end{frontmatter}

\section{Introduction}

We consider the infinite-source Poisson process with random transmission rate
defined by
\begin{equation} \label{eq:defgginfini}
  \Xgen(t) = \sum_{\ell \in \Nset} U_\ell  \ind_{\{ t_\ell \leq t < t_\ell + \eta_\ell\}},\qquad t\geq 0,
\end{equation}
where the arrival times $\{t_\ell\}_{\ell\geq0}$ are the points of a unit-rate
  homogeneous Poisson process on the positive half-line, independent of the
  initial conditions; and
the durations and transmission rates $\{(\eta_\ell,U_\ell)\}$ are independent and identically distributed
  random variables with values in $(0,\infty)\times\Rset$ and independent of the
  Poisson process and of the initial conditions.
This process was considered by Resnick and Rootz\'en \cite{resnickrootzen2000} and
Mikosch \textit{et al.} \cite{mikoschresnickrootzenstegeman2002}, among others. The $M/G/\infty$
queue is a special case, for $U_\ell \equiv 1$. An important motivation for
the infinite-source Poisson process is to model the instantaneous rate of the
workload going though an internet link. Although overly simple models are
generally not relevant for internet traffic at the packet level, it is
generally admitted that rather simple models can be used for higher-level (the
so-called \textit{flow level}) traffic such as TCP or HTTP sessions, one of
them being the infinite-source Poisson process (see Barakat \textit{et al.} \cite{flowlevel2002}).  One
way to empirically analyse internet traffic at the flow level using the
infinite-source Poisson process would consist in retrieving all the variables
$\{t_\ell,\eta_\ell,U_\ell\}$ involved in the observed traffic during a given period of
time, but this would require the collection of all the relevant information in the
packets headers (such as source and destination addresses) for the purpose of separating the
aggregated workload into transmission rates at a pertinent level;
see Duffield \textit{et al.}~\cite{flowlevelstat2002} for many insights into this
problem.

It is well known that heavy tails in the durations $\{\eta_k\}$ result in
long-range dependence of the process $\Xgen(t)$.
Long-range dependence can be defined by the regular variation of the
autocovariance of the process or more generally by the regular variation of the
variance of the integrated process:
\[
  \var \biggl( \int_0^t \Xgen(s)  \d s \biggr) = L(t) t^{2H},
\]
where $L$ is a slowly varying function at infinity and $H>1/2$ is often refered
to as the \textit{Hurst index} of the process. For the infinite-source Poisson
process, the Hurst index $H$ is related to the tail index $\alpha$ of the
durations by the relation $H = (3-\alpha)/2$.  The long-range dependence
property has motivated many empirical studies of internet traffic and
theoretical ones concerning its impact on queuing (these questions are studied
in the $M/G/\infty$ case in Parulekar and Makowki \cite{parulekarmakowki1997}).

However, to the best of our knowledge, no statistical procedure to
estimate $H$ has been rigorously justified. It is the aim of this
paper to propose an estimator of the Hurst index of the
infinite-source Poisson process, and to derive its statistical properties.  We
propose to estimate $H$ (or equivalently $\alpha$) from a path of the
process $\Xgen(t)$ over a finite interval $[0,T]$, observed either
continuously or discretely. In practice this can be done by counting
all the packets going through some point of the network and then
collecting local traffic rate measurements.  Our estimator is based on
the so-called wavelet coefficients of a path. There is a wide
literature on this methodology for estimating long-range dependence,
starting as long ago as Wornell and Oppenheim \cite{wornelloppenheim1992}, but
we are not aware of rigorous results for non-Gaussian or non-stable
processes. The main contribution of this paper is thus the proof of
the consistency of our estimator. We also investigate the rate of
convergence of the estimator in the case $\alpha>1$. If the process is
observed continuously, the rate of convergence is good. In the case of
discrete observations, the rate is much smaller. Also, the choice of
the tuning parameters of the estimators is much more restricted in the
latter case, and practitioners should perhaps be aware of this; see
Section~\ref{sec:rate-conv-stable} for details.

The process $\Xgen$ is formally defined in Section~\ref{sec:basic-prop-model}. We
state our assumptions and, using a point-process representation of $\Xgen$, we
establish some of its main properties. The wavelet coefficients are defined
and the scaling property of their variances is obtained in
Section~\ref{sec:wavelet-coefficients}.  The estimator is defined and its
properties are established in Section~\ref{sec:estimation}. The \hyperref[sec:technical-results]{Appendix}
contains technical lemmas.

\section{Basic properties of the model}
\label{sec:basic-prop-model}

\subsection{Assumptions}

We now introduce the complete assumption on the joint distribution of the
transmissions rates and durations.
\begin{hypo} \label{hypo:modele}
\textup{(i)} The random vectors $\{(\eta,U),\  (\eta_\ell,U_\ell),\
  \ell\in\Zset\}$ are independent  with common distribution $\nu$ on $(0,\infty)\times\Rset$ and
  independent of the homogeneous Poisson point process on the real line with
  points $\{t_\ell\}_{\ell\in\Zset}$ such that $t_{\ell}<t_{\ell+1}$ for all $\ell$ and
  $t_{-1}<0\leq t_0$.

\textup{(ii)} There exists a positive integer $p^*$ such
  that $\esp[|U|^{p^*}] <\infty$.

\textup{(iii)} There exist a real number
  $\alpha \in (0,2)$ and  positive
   functions $L_0,\ldots,L_{p^*}$ slowly varying at infinity such that, for all $t>0$ and  $p=0,\dots,p^*$,
  \begin{equation}
    \label{eq:tailbehavior}
    H_p(t) := \esp \bigl[ |U|^p \ind_{\{\eta >t\}} \bigr] = L_p(t) t^{-\alpha}  .
  \end{equation}
\end{hypo}

Since $\eta>0$, the functions $H_p$ are continuous at zero and $H_p(0)
= \esp[|U|^p]$.  Condition~(\ref{eq:tailbehavior}) is equivalent to
saying that the functions $H_p,p=0,1,\dots,p^*,$ are regularly
varying with index $-\alpha$.  If $\alpha>1$ and $p^* \geq2$,
Assumption~\ref{hypo:modele} and Karamata's theorem imply the
following asymptotic equivalence:
\begin{equation}  \label{eq:momentueta3}
  \esp [U^2 \{\eta - t\}_+ ] =
\esp \biggl[ U^2 \int_{v=t}^\infty \ind_{\{v < \eta\}}  \d v \biggr]
= \int_{v=t}^\infty  H_2(v) \d v   \sim \frac 1 {\alpha-1} L_2(t) t^{1-\alpha}  .
\end{equation}

\begin{rem}
  Assumption~\ref{hypo:modele} will be used with $p^*=2$ to prove the
  regular variation of the autocovariance function of the process
  $\Xgen$ and with $p^*=4$ to prove consistency of our estimators. It
  can be related to the theory of multivariate regular variation
 (see, for instance, Maulik \textit{et al.} \cite{maulikresnickrootzen2002}). But the
definitions of multivariate regular variation involve vague
convergence and do not necessarily ensure the convergence of moments
required here.
\end{rem}

\begin{rem}
  We do not assume that $U$ is non-negative. This allows us to consider
  applications other than teletraffic modelling. For instance, the process $\Xgen$
  could be used to model the volatiliy of some financial time series.
\end{rem}

\begin{rem}
We will often have to separate the cases $\esp[\eta]=\infty$
and $\esp[\eta]<\infty$.  These cases are respectively implied by $\alpha<1$
and $\alpha>1$. If $\alpha=1$, the finiteness of $\esp[\eta]$ depends on the
precise behaviour of $L_0$ at infinity.
\end{rem}

\begin{xmpl}
Assumption~\ref{hypo:modele} implies in particular that the tail of the
distribution of $\eta$ is regularly varying with index $\alpha$. This in
turns implies Assumption~\ref{hypo:modele} if $U$ and $\eta$ are independent
and $\esp[ |U|^{p^*}] < \infty$, in which case the functions $L_p$ differ by
a multiplicative constant.
\end{xmpl}

\begin{xmpl}
  Assumption~\ref{hypo:modele} also holds in the following case which is  of
  interest in teletraffic modelling.  In a TCP/IP traffic context, $\eta$ and
  $U$ represent respectively the duration of a download session and its
  intensity (bit rate). Then $W := U\eta $ represents the amount of transmitted
  data. We assume that, for some $u_0>0$, there exist two regimes, $U\geq u_0$ (xDSL/LAN/cable
  connection) and $U\in(0,u_0)$ (RTC connection), such that the following statements hold:
\begin{itemize}
\item The distribution of $W$ given $U=u \geq u_0$ is heavy-tailed and
  independent of $u$: $\pr(W \geq w | U = u)= L(w) w^{-\alpha}$.
\item The distribution of $W$ given $U=u\in(0,u_0)$ is light-tailed
  uniformly with respect to $u$. For instance, we assume expontially decaying
  tails, $\pr(W \geq w | U = u) \leq \exp(-\beta w^{-\gamma})$, for some $\beta>0$ and $\gamma>0$.
\end{itemize}
An explicit example for two such regimes is obtained when the conditional density of $W$ given $U=u$ is equal to
$\alpha w^{-\alpha-1}\ind_{\{w\geq1\}}$ if $u\geq u_0$ and $\exp(-w)$ if $u< u_0$.

Concerning the distribution of $U$ we only assume that:
\begin{itemize}
\item  $\pr(U\geq u_0)>0$, $\esp[|U|^{-\alpha-\epsilon}] < \infty$ for some
$\epsilon>0$, and $\esp[|U|^{p^*}] < \infty$.
\end{itemize}
Then (\ref{eq:tailbehavior}) holds for $p\leq p^*$. Indeed,
\begin{eqnarray}\label{eq:highrateEq}
\nonumber
   \esp\bigl[U^p\ind_{\{\eta > t\}} \ind_{\{U \geq u_0\}}\bigr]
   &=& \esp\bigl[U^p\ind_{\{W > Ut\}} \ind_{\{U \geq u_0\}}\bigr] \\
\nonumber
  & = & \esp \bigl[U^p L(Ut)(Ut)^{-\alpha} \ind_{\{U\geq u_0\}}\bigr]  \\
  & = & L(t) t^{-\alpha} \esp \bigl[U^{p- \alpha} L(Ut)/L(t) \ind_{\{U\geq u_0\}}\bigr]  .
\end{eqnarray}
Since $L$ is slowly varying at infinity, $\lim_{t\to \infty} L(ut)/L(t) = 1$,
uniformly with respect to $u$ in compact sets of $(0,+\infty)$, and there exists
 $t_0>0$ such that, for $u\geq u_0$, $t \geq t_0$,
\[
  \frac{L(ut)}{L(t)} \leq (1+\alpha) u^{\alpha/2}  ;
\]
see, for example, Resnick (\cite{resnick1987}, Proposition 0.8).
Then, by the dominated convergence theorem,
\begin{equation}
\label{eq:highrate}
  \lim_{t \to \infty} \esp \bigl[U^{p- \alpha}  L(Ut)/L(t) \ind_{\{U\geq u_0\}}\bigr] =
  \esp\bigl[U^{p-\alpha} \ind_{\{U\geq u_0\}}\bigr].
\end{equation}
Consider now the low-bit-rate regime. Since, for all $x>0$, $\exp\{-\beta x^\gamma\}\leq C x^{-\alpha-\epsilon}$ for
some positive constant $C$, we have
\[
\esp\bigl[U^p\ind_{\{\eta > t\}}\ind_{\{U< u_0\}}\bigr] \leq \esp \bigl[U^p \exp\{-\beta (Ut)^\gamma\}\ind_{\{U< u_0\}}\bigr]
\leq C t^{-\alpha-\epsilon} \esp \bigl[U^{p - \alpha - \epsilon} \ind_{\{U< u_0\}} \bigr] .
\]
Using the assumption on $U$, since $p\geq0$, the rightmost expectation
in the previous display is finite and we obtain that
\[
  \lim_{t \to \infty} t^{\alpha} L^{-1}(t) \esp\bigl[U^p\ind_{\{\eta > t\}} \ind_{\{U< u_0\}} \bigr] = 0 .
\]
Together with (\ref{eq:highrateEq}) and (\ref{eq:highrate}), this implies that, as $t\to\infty$,
$\esp[U^p\ind_{\{\eta > t\}}]t^\alpha \sim L(t) \esp[U^{p-\alpha}\times \ind_{\{U\geq u_0\}}]$
hence is slowly varying and Assumption~\ref{hypo:modele} holds.
\end{xmpl}

\subsection{Point-process representation and stationary version}
\label{sec:ppp}

Let $\mcn$ denote a Poisson point process on a set $E$ endowed with a
$\sigma$-field $\mce$ with intensity measure $\mu$, that is, a random
measure such that for any disjoint $A_1,\ldots,A_p$ in $\mce$,
$\mcn(A_1),\dots,\mcn(A_p)$ are independent random variables with
Poisson law with respective parameters $\mu(A_i)$, $i=1,\dots,p$.  The
main property of Poisson point processes that we will use is the
following cumulant formula (see, for
instance, Resnick \cite{resnick1987}: Chapter~3). For any positive integer $p$ and
functions $f_1,\dots,f_p$ such that $\int |f_i| \d \mu<\infty$ and
$\int |f_i|^p \d \mu < \infty$ for all $i=1,\dots,p$, the $p$th-order
joint cumulant of $\mcn(f_1),\dots,\mcn(f_p)$ exists and is given by
\begin{equation} \label{eq:cum}
   \cum(\mcn(f_1),\dots,\mcn(f_p) ) = \int f_1 \cdots f_p  \d \mu  .
\end{equation}

Let $\NS$ be the point processes on $\Rset\times(0,\infty)\times\Rset$ with points
$(t_\ell,\eta_\ell,U_\ell)_{\ell\in\Zset}$, that is $\NS = \sum_{\ell\in\Zset} \delta_{t_\ell,\eta_\ell,U_\ell}$.
Under Assumption~\ref{hypo:modele}(i), it is a Poisson point
process with intensity measure  $\Leb \otimes \nu$, where $\Leb$
is the Lebesgue measure on $\Rset$. For
$t,u \in \Rset$, define
\begin{eqnarray*}
  A_t & =& \{(s,v) \in \Rset\times\Rset_+ \mid s \leq t < s + v \}  ,\\
  B_u & =& \{\lambda u\mid\lambda\in[1,\infty)\} .
\end{eqnarray*}

We can now show that if $\esp[\eta]<\infty$, then one can define a
stationary version for $\Xgen$ and provide its second-order properties.
\begin{prop}\label{prop:XsAllPropoerties}
  If Assumption~\textup{\ref{hypo:modele}(i)} holds and
  $\esp[\eta]<\infty$, then the process
\begin{equation} \label{eq:defgginfinistationnaire}
  \XS(t) = \sum_{\ell \in \Zset} U_\ell \ind_{\{ t_\ell \leq t < t_\ell + \eta_\ell\}}
\end{equation}
is well defined and strictly stationary.  It has the point-process
representation
\begin{equation}\label{eq:XSPPrep}
  \XS(t) = \int_{0}^\infty \NS(A_{t}\times B_{u})  \d u - \int_{-\infty}^0
  \NS(A_{t}\times B_{u})  \d u  .
\end{equation}

Let $K_0 = \sup\{ \ell > 0 \mid t_{-\ell}+ \eta_{-\ell} >0\}$, $\tilde U_\ell = U_{-\ell}$ and
$\tilde \eta_\ell = \eta_{-\ell}+t_{-\ell}$. Then, for all $t\geq0$,
\begin{equation}\label{eq:XSisXplusInitCond}
\XS(t)=\sum_{\ell=1}^{K_0} \tilde{U}_\ell \ind_{\{ t < \tilde \eta_\ell\}} + \Xgen(t).
\end{equation}

If, moreover, $p^*\geq2$, then $\XS$ has finite
variance and
\begin{eqnarray*}
  \esp[\XS(t)] & = & \esp[U \eta ]  , \\
  \cov(\XS(0), \XS(t)) & = & \esp[ U^2 (\eta - t)_+] = \int_t^\infty H_2(v) \d v .
\end{eqnarray*}
\end{prop}

\begin{rem}
  Note that  if $\alpha>1$, then $\esp[\eta] < \infty$ and, by Karamata's
  theorem,
\[
  \cov(\XS(0), \XS(t)) \sim \frac1{\alpha-1} L_2(t) t^{1-\alpha}  , \qquad {t
    \to +\infty}  .
\]
\end{rem}

\begin{pf}
  The number of non-vanishing terms in the
  sum~(\ref{eq:defgginfinistationnaire}) is $\NS(A_t\times \Rset)$ and has a
  Poisson distribution with mean $\esp \int_{\Rset}\ind_{A_t}(s,\eta)  \d s =
  \esp[\eta]$.  Thus $\XS$ is well defined and stationary since $\NS$ is
  stationary. The number of indices $\ell>0$ such that
  $t_{-\ell}+\eta_{-\ell}>0$ is $N_S(A_0\times \Rset)$, hence if $K_0$ is the
  largest of those $\ell$s, it is almost surely finite and
\[
 \sum_{t_\ell<0} U_\ell \ind_{\{t_\ell \leq t < t_\ell + \eta_\ell\}}  =
  \sum_{\ell=1}^{K_0} \tilde{U}_\ell \ind_{\{ t < \tilde \eta_\ell\}}  .
\]
Hence~(\ref{eq:XSisXplusInitCond}).

The point-process representation~(\ref{eq:XSPPrep}) and
formule~(\ref{eq:cum}) and~(\ref{eq:momentueta3}) finally yield the
given expressions for the mean and covariance.
\end{pf}

Relation~(\ref{eq:XSisXplusInitCond}) shows that the stationary
version $\XS$ can be defined by changing the initial condition of the
system. More generally, one could consider \textit{any} initial
conditions, that is, any process defined as on the right-hand side
of~(\ref{eq:XSisXplusInitCond}) with $K_0$ and $\tilde \eta_\ell,
\ell>0$ finite.  Since the initial conditions almost surely vanish
after a finite period, they have a negligible impact on the estimation
procedure.  Thus, our result on $\Xgen$ easily generalizes to any such
initial conditions, and, in particular, to the stationary version
$\XS$, when it exists.

Applying similar arguments as those used for showing Proposition~\ref{prop:XsAllPropoerties}, we
obtain:
\begin{prop} \label{prop:secondorder}
The process $\Xgen$ admits a point-process representation
\begin{equation}\label{eq:XZPPrep}
  \Xgen(t) = \int_{0}^\infty \NS(A_{t}^+\times B_u)  \d u - \int_{-\infty}^0
  \NS(A_{t}^+\times B_u)  \d u  ,
\end{equation}
where $A_t^+=A_t\cap \Rset_+^2$.

  If Assumption~\ref{hypo:modele} holds with $p^*\geq 2$, then the process $\Xgen$ is
  non-stationary with expectation and autocovariance function given, for $s \leq
  t$, by
\begin{eqnarray*}
  \esp[\Xgen(t)] &=& \esp[U (\eta \wedge t)]  , \\
\cov(\Xgen(s), \Xgen(t)) &=& \esp[ U^2 \{s - (t-\eta)_+\}_+] =  \int_{t-s}^t H_2(v) \d v .
\end{eqnarray*}
\end{prop}

By the uniform convergence theorem for slowly varying functions, the
following asymptotic equivalence of the covariance holds.  For any
$\alpha\in(0,2)$ and all $t>s>0$, as $T\to\infty$,
\begin{equation}
\label{eq:asympcovnonstat}
 \cov(\Xgen(Tt), \Xgen(Ts)) \sim C L_2(T)  T^{1-\alpha}
\end{equation}
with $C= \int_{t-s}^s v^{-\alpha} \d v$.

In accordance with the notation in use in the context of long-memory
processes, we can define the Hurst index of the process $\Xgen$ as
$H=(3-\alpha)/2$, because the variance of the process integrated
between 0 and $T$ increases as $T^{2H}$. If $\alpha<1$, then $H>1$.
This case has been considered, for instance, by
Resnick and Rootz\'en \cite{resnickrootzen2000}.

\section{Wavelet coefficients}
\label{sec:wavelet-coefficients}

\subsection{Continuous observation}
\label{sec:defin-stat-vers}

Let $\psi$ be a bounded real-valued function with compact support in $[0,M]$
and such that
\begin{equation} \label{eq:momentnul}
  \int_{0}^M \psi(s)  \d s = 0  .
\end{equation}
For integers $j\geq0$ and $k\in\Zset$, define
\begin{equation}\label{eq:psijk}
  \psi_{j,k}(s) = 2^{-j/2} \psi(2^{-j}s-k).
\end{equation}
The wavelet coefficients of the path are defined as
\begin{equation} \label{eq:defbdjk}
d_{j,k} = \int_0^{\infty}
  \psi_{j,k}(s) \Xgen(s)  \d s
\end{equation}
(see, for example, Cohen \cite{cohen2003}).
Assume that a path of the process $\Xgen$ is observed continuously between times
0 and $T$. Since $\psi_{j,k}$ has support in $[k2^j,(k+M)2^j]$, the
coefficients $d_{j,k}$ can be computed for all $(j,k)$ such
that $T2^{-j}\geq M$ and $k=0,1,\dots,T2^{-j}-M$.

According to Lemma~\ref{lem:unobs}, one may define, for all $j$ and $k$,
\begin{equation} \label{eq:defunobservable}
  d_{j,k}^S = \sum_{\ell\in\Zset} U_\ell \int_{t_\ell}^{t_\ell+\eta_l}
  \psi_{j,k}(s)\d s  .
\end{equation}
As stated in Lemma~\ref{lem:unobs}, if $\esp[\eta] <\infty$, we have
$d_{j,k}^S=  \int_0^{\infty} \psi_{j,k}(s) \XS(s) \d s$.
Nevertheless, even if $\esp[\eta] = \infty$, the sequence of coefficients at a
given scale $j$, $\{d^S_{j,k},k\in\Zset\}$, is stationary.
Moreover, the definition~(\ref{eq:defunobservable}) yields:
\begin{lem} \label{lem:varwavcoefstationnaire}
  Let Assumption~\ref{hypo:modele} hold with $p^*\geq2$. We have
   \begin{equation}  \label{eq:varwavcoeffstationaire}
     \esp[d_{j,k}^S] = 0  , \qquad \var(d_{j,k}^S) = \mathcal{L} (2^j)
      2^{(2-\alpha)j}  ,
 \end{equation}
where
 \begin{equation} \label{eq:defmathcall}
   \mathcal L (z) := z^{\alpha} \int_0^\infty \int_{-\infty}^\infty
   \biggl( \int_{-\infty}^\infty \biggl\{ \int_t^{t+v z^{-1}} \psi(s)
       \d s \biggr\} ^2  \d t \biggr)  w^2 \nu(\d v, \d w)
 \end{equation}
 is slowly varying as $z\to\infty$. More precisely, we have
 the asymptotic equivalence
 \begin{equation} \label{eq:mathcallequiv}
\mathcal L(z) \sim C_\mcl L_2(z) \qquad\mbox{as }z\to\infty ,
\end{equation}
with $C_\mcl = \alpha\int_0^\infty \int_{-\infty}^\infty \{ \int_x^{x+y} \psi(s)  \d s
 \}^2  \d x \, t^{-\alpha-1} \d t >0$.
\end{lem}

The proof of~(\ref{eq:varwavcoeffstationaire}) is a
straightforward application of (\ref{eq:cum}), and the proof of
the asymptotic equivalence~(\ref{eq:mathcallequiv}) is obtained by
standard arguments on slowly varying functions. A detailed proof can be
found in Fa\"y \textit{et al.} \cite{fayroueffsoulier2005arxiv}.

\subsection{Wavelet coefficients in discrete time}\label{sec:discr-wavel-coeff}

Let $\phi$ be a bounded $\Rset\to\Rset$ function with compact support
included in $[-M+1,1]$ and such that
\begin{equation}\label{eq:phiInt}
\sum_{k\in\Zset} \phi(t-k) = 1,\qquad t\in\Rset.
\end{equation}
Let $\interp$ denote the operator defined on the set of functions
$x\dvtx\Rset\to\Rset$ by
\begin{equation}\label{eq:interpStep}
\interp[x](t)=\sum_{k\in\Zset}x(k) \phi(t-k).
\end{equation}
The  wavelet coefficients of $x$ are then defined as the wavelet coefficients of
$\interp[x]$.

From a computational point of view, it is convenient to chose $\phi$
and $\psi$ to be the so-called father and mother wavelets of a multiresolution
analysis; see, for instance, Meyer~\cite{meyer1992}. The simplest choice
is to take $\phi$ and $\psi$ to be associated with the Haar system, in which case
$M=1$, $\phi=\ind_{[0,1)}$ and $\psi=\ind_{[0,1/2)}-\ind_{[1/2,1)}$.

If the process $\Xgen$ is observed discretely, we denote its wavelet coefficients
by
\begin{equation}\label{eq:djk_discrete_case*}
  d^D_{j,k} = \int \psi_{j,k}(s) \interp[\Xgen](s)  \d s  .
\end{equation}
If we observe
$\Xgen(0),\Xgen(1),\dots,\Xgen(T-1)$, for some positive integer $T$,
we can compute $d^D_{j,k}$ for all $j,k$ such that $0 \leq k \leq 2^{-j}
(T-M+1) - M$.
Roughly, for $2^j \geq T/M$, no coefficients can be computed and if $2^j < T/M$
the number of computable wavelet coefficients at scale $2^{-j}$ is of
order $T2^{-j}+1-M$ for $j$ and $T$ large.


\begin{rem}
  Observe that the choice of time units is unimportant here.  Indeed, in
  Assumption~\ref{hypo:modele}, changing the time units simply amounts to
  adapting the slowly varying functions $L_k$ and the rate of the arrival process
  $\{t_k\}$. Clearly these adaptations do not modify our results since
  precise multiplicative constants are not considered.
\end{rem}

\subsection{Averaged  observations}

We describe now a third observation scheme for which our results can easily be
extended.  Suppose that $T$ is a positive integer and that we observe local
averages of the trajectory
\[
\overline{\Xgen}(k):= \int_k^{k+1} \Xgen(t)  \d t = \int \Xgen(t)  \phi_H(t-k)
 \d t,\qquad k=0,1,\dots, T-1  ,
\]
where $\phi_H:=\ind_{[0,1]}$ is the Haar wavelet.  Let $\interpbis$ denote
the operator on locally integrable functions $x$ defined by
\[
\interpbis[x](t)=\sum_{k\in\Zset}  \biggl(\int x(s) \phi_H(s-k) \d s \biggr)
\phi(t-k).
\]
For this observation scheme, as in Section~\ref{sec:discr-wavel-coeff}, one
may compute the wavelet coefficients of the function $\interpbis[\Xgen]$ at all
scale and location indices $(j,k)$ such that $0 \leq k \leq 2^{-j} (T-M+1) -
M$.  If $\phi=\phi_H$
and  $\psi$ is the Haar mother wavelet,
$\psi=\ind_{[0,1/2)}-\ind_{[1/2,1)}$, then
the wavelet coefficients of $\interpbis[\Xgen]$ are precisely the continuous
wavelet coefficients defined in~(\ref{eq:defbdjk}).  For any other choice of
$\phi$ and $\psi$, this is no longer true. We will not treat this case, but all
our results can be extended at the cost of further technicalities.

\section{Estimation}
\label{sec:estimation}

Tail index estimation methods do not seem appropriate here for estimating the
parameter $\alpha$. Indeed, $\alpha$ is the tail index of the
unobserved durations $\{\eta_k\}$, whereas the observed process $\Xgen(t)$ always
has finite variance ($\esp [|\Xgen(t)|^p] < \infty$ if and only if
$\esp[|U^p|]<\infty$ and the marginal distribution of $\Xgen(t)$ is Poisson if
$U = 1$ almost surely).  But as shown by Proposition~\ref{prop:secondorder}, $\alpha$
is related to the second-order properties of the process: the coefficient $H =
(3-\alpha)/2$ can be viewed as its Hurst index, that is, $H$ governs the rate of
decay of the autocovariance function of the process.  Therefore it seems
natural to use an estimator of the Hurst index.

\subsection{The estimator}

Lemma~\ref{lem:varwavcoefstationnaire} provides the rationale for a
minimum contrast estimator of $\alpha$ which is related to the local Whittle
estimator; cf. K\"unsch \cite{kuensch1987}.
Let $\bd_{j,k}$ denote the wavelet coefficients which are actually
  available; these may be obtained from continuous-time
  ($\bd_{j,k}=d_{j,k}$) or discrete-time ($\bd_{j,k}=d^D_{j,k}$)
  observations.
Let $\Delta$ be a set of indices $(j,k)$ of available wavelet
coefficients.
Denote the mean scale index over $\Delta$ by
\[
  \delta := \frac{1}{\# \Delta} \sum_{(j,k) \in \Delta} j  .
\]
The reduced local Whittle contrast function is
\begin{equation}\label{eq:whittlecontrast}
 \hat W(\alpha') =  \log \Biggl( \sum_{(j,k)\in\Delta}
      \frac{\bd^{2}_{j,k}}{2^{(2-\alpha')j}} \Biggr) + \delta \log(2)
    (2-\alpha')  .
\end{equation}
The local Whittle estimator of $\alpha$ is then defined as
\begin{equation}\label{eq:defestimator}
  \hat \alpha := \arg\min_{\alpha'\in (0,2)}  \hat W(\alpha').
\end{equation}
In order to simplify the proof of our result, we henceforth take $\Delta$ to be of
the form
\[
  \Delta  = \{(j,k)  ; J_0 < j \leq J_1  ,  0 \leq k \leq \nj_j-1 \}  ,
\]
with
$ J  =  \max\{j ; 2^j \leq (T-M+1)/(M+1) \}$,
$ n_j  = 2^{J-j} $
and integers $J_0$ and $J_1$  such that
\begin{equation}\label{eq:basicCondonJs}
  0<J_0 < J_1 \leq J.
\end{equation}
The sequence of integers $J$ depends on $T$ in such a way that $2^J \asymp T$.
Note that the dependence of the sequences $J$, $J_0$, $J_1$, $n_j$
etc. on $T$ is suppressed in our notation.

\subsection{Consistency}

Our estimator is consistent in the potentially unstable case, that is when
$\alpha$ is not assumed to be in $(1,2)$, provided that the assumptions on the
functions $\phi$ and $\psi$ are strengthened. We assume that
\begin{equation} \label{eq:2ndmomentnul}
  \int_{-\infty}^\infty s\psi(s)  \d s = 0  ,
\end{equation}
and there exist constants $a$ and $b$ such that, for all $t\in\Rset$,
\begin{equation}
\label{eq:phiIntLin}
\sum_{k\in\Zset} k \phi(t-k) = a + bt  .
\end{equation}
These conditions are not satisfied by the Haar wavelet, but hold for any
Daubechies wavelets; see Cohen \cite{cohen2003}.

\begin{thm} \label{theo:stable}
  Let Assumption \ref{hypo:modele} hold with $p^*\geq4$.  Assume that $J_0$ and $J_1$
  depend on $T$ in such a way that
  \begin{eqnarray}
 \lim_{T\to\infty} J_0 & = &    \lim_{T\to\infty} (J_1-J_0) = \infty  , \label{eq:rateJ0J1} \\
 \limsup_{T\to\infty} {J_0}/J &<& 1/\alpha ,\label{eq:rateJ0J} \\
 \limsup_{T\to\infty} J_1/J &<& 1/(2-\alpha)  . \label{eq:rateJ1J}
  \end{eqnarray}
  Then $\hat \alpha$ is a consistent estimator of $\alpha$. Moreover if
  $\alpha\in(1,2)$, then
  conditions~(\ref{eq:2ndmomentnul}),~(\ref{eq:phiIntLin})
  and~(\ref{eq:rateJ1J}) are not necessary for the same result to hold.
\end{thm}


\begin{rem}\label{rem:universalEstimator}
  Conditions (\ref{eq:rateJ0J1}), (\ref{eq:rateJ0J})
  and (\ref{eq:rateJ1J}) are satisfied by the choice $J_0 = \lfloor J/2
  \rfloor$ and $J_1 = \lfloor J/2+\log(J) \rfloor$.
\end{rem}

\begin{pf*}{Proof of Theorem~\ref{theo:stable}}
  For clarity of notation, we denote $\sum_j = \sum_{j=J_0+1}^{J_1}$, $\Delta_j
  := \{k\dvtx (j,k)\in\Delta\}$ and $\# \Delta_j = n_j$.  Elementary computations
  give
\begin{equation}\label{eq:delta}
\delta = J_0 + 2 + (J_0-J_1)/(2^{J_1-J_0}-1)
\end{equation}
so that $\delta-(J_0 + 2)\to 0$ under (\ref{eq:rateJ0J1}).  By Karamata's
representation theorem, the slowly varying function $\mcl$ defined
in (\ref{eq:defmathcall}) can be written as
\[ \label{sec:representationmathcall}
  \mcl(z) = c  \bigl(1 + r(z)\bigr)  \exp \biggl\{\int_1^z \frac{\ell(s)}{s} \d
    s\biggr\}  ,
\]
with $c>0$ and $\lim_{z \to \infty} \ell(z) = \lim_{z \to \infty} r(z) = 0$.
Define $\mcl_0(z) = c \exp \{\int_1^z s^{-1} \ell(s) \d s\} $,
$r^*(z) = \sup_{z' \geq z} |r(z')|$ and $\ell^*(z) = \sup_{z' \geq z}
|\ell(z')|$. The functions $r^*$ and $\ell^*$ are non-increasing and tend to
zero at infinity.  We now introduce some notation that will be used throughout the
proof:
  \begin{eqnarray*}
    W(\alpha') &\hspace*{2.5pt}=& \log \biggl( \sum_j 2^{(\alpha'-\alpha)j} n_j
      \mcl(2^j) \biggr) + \delta \log(2) (2-\alpha')      ,  \\
W_0(\alpha') &\hspace*{2.5pt}=& \log \biggl( \sum_j 2^{(\alpha'-\alpha)j} n_j
      \mcl_0(2^j) \biggr) + \delta \log(2) (2-\alpha')
     ,  \\
  w_{j,0}(\alpha') &:=& \frac{2^{(\alpha'-\alpha)j} n_j \mathcal L_0(2^j)}
    {\sum_{j'} 2^{(\alpha'-\alpha)j'} n_{j'} \mathcal L_0(2^{j'}) }  ,
    \qquad
    w_{j}(\alpha') := \frac{2^{(\alpha'-\alpha)j} n_j \mathcal L(2^j)}
    {\sum_{j'} 2^{(\alpha'-\alpha)j'} n_{j'} \mathcal L(2^{j'}) }  , \\
   v_j &\hspace*{2.5pt}=& \mathcal L(2^j)  2^{(2-\alpha)j}  , \qquad \Lambda_j = \nj_j^{-1}
    \sum_{k=0}^{\nj_j-1} \{     v_j^{-1} (\bd_{j,k})^2 - 1 \}  , \\
    E(\alpha') &:=& \sum_j w_j(\alpha') \Lambda_j  .
  \end{eqnarray*}
We have
\[
W(\alpha')-W_0(\alpha')=\log\biggl( 1+\frac{\sum_j
2^{(\alpha'-\alpha)j} n_j \mcl_0(2^j)r(2^j)}{\sum_j
2^{(\alpha'-\alpha)j} n_j \mcl_0(2^j)} \biggr)  .
\]
Here the fraction inside the logarithm is bounded by
$r^*(2^{J_0})$, thus, for $J$ large enough,
\[
\sup_{\alpha'}|W(\alpha')-W_0(\alpha')| \leq C r^*(2^{J_0}).
\]
Standard algebra yields
\begin{eqnarray*}
W_0'(\alpha') &=& \log 2 \sum_j w_{j,0}(\alpha') (j-\delta) \\
W_0''(\alpha') &=& \log^2 (2) \sum_j w_{j,0}(\alpha') \biggl( j -
  \sum_{j'} w_{j',0}(\alpha') j' \biggr)^2 .
\end{eqnarray*}
By Lemma~\ref{lem:svfunc}, under~(\ref{eq:rateJ0J1}),
\[
  \lim_{T\to\infty} W'_0(\alpha) = 0  ,\qquad  \lim_{T\to\infty} W''_0(\alpha) =
  2   .
\]
Thus, there exist $\eta>0$ and $\zeta>0$ such that
\[
  \liminf_{T\to\infty} \inf_{\alpha' \in (\alpha-\eta,\alpha+\eta)}
  W''_0(\alpha') > \zeta  .
\]
This implies that, for large $T$ and some positive constant $c$,
\begin{equation} \label{eq:minoration}
  W(\hat \alpha) - W(\alpha) \geq W_0'(\alpha) \log(2) (\hat\alpha-\alpha) +
  c(\hat \alpha - \alpha)^2 - 2r^*(2^{J_0}).
\end{equation}
Since $W'_0(\alpha)\to0$ and $|\hat\alpha-\alpha|\leq 2$, this implies that, for
all $\epsilon>0$,
\begin{equation} \label{eq:sufficientconsistency}
  \limsup_{T\to\infty} \pr\bigl( (\hat \alpha - \alpha)^2 > \epsilon\bigr)
  \leq \limsup_{T\to\infty} \pr\bigl( W(\hat \alpha) - W(\alpha) \geq
  c\epsilon\bigr)  .
  \end{equation}
  Write
\begin{eqnarray}\label{eq:majo}
   \hat W(\alpha') & = & W(\alpha') + \log \{ 1 + E(\alpha') \}  ,
 \nonumber \\
    W(\hat \alpha) - W(\alpha) & =& \hat W(\hat \alpha) - \hat W(\alpha)
   - \log \{ 1 + E(\hat \alpha) \}
   + \log \{ 1 + E(\alpha) \}  \nonumber \\
   & \leq& 2 \sup_{\alpha' \in (0,2)} | \log \{ 1 +
       E(\alpha') \} |  .
 \end{eqnarray}
 Consistency will follow from (\ref{eq:sufficientconsistency})
 and~(\ref{eq:majo}) provided that we can prove that\break $\sup_{\alpha' \in (0,2)}
 |E(\alpha')| =o_P(1)$.  If $\alpha>1$, take $\varepsilon \in (0,
 (\alpha-1)/2)$ such that $\limsup J_0/J < 1/(\alpha + \varepsilon)$, which is
 possible by assumption~(\ref{eq:rateJ0J}).  Define
\begin{equation}  \label{eq:conditionJ11}
  J_2 =
  \cases{
    J_1,  &\quad if $ \alpha  \leq 1   $, \cr
    J_1 \wedge [J/(\alpha+\varepsilon)], &\quad  if $ \alpha >1  $,
  }
\end{equation}
so that, for $T$ large enough, $J_0< J_2 \leq J_1$.
Write
\[
  E(\alpha')  = \sum_{j = J_0 + 1}^{J_2} w_j(\alpha') \Lambda_j + \sum_{j
    = J_2 + 1}^{J_1} w_j(\alpha') \Lambda_j =: E_1(\alpha') + E_2(\alpha')  ,
\]
with the convention that $\sum_{j = J_2 + 1}^{J_1}=0$ if $J_2=J_1$.
By Lemma~\ref{lem:blob},
\begin{equation}
\sup_{\alpha'\in(0,2)} |E_1(\alpha')|
= O_P(2^{-\xi_1 J}), \label{eq:borneE1}
\end{equation}
for some positive $\xi_1$.
Now treat $E_2$ for $\alpha > 1$ and $J_2=[J/(\alpha+\varepsilon)] > J_1$. For all $\alpha'
\in (0,2)$, we have $\alpha' - \alpha - 1 < -2\varepsilon$.  Since $\mathcal L$ is
slowly varying, we obtain, for some positive constant $C$, for all
$j=J_2+1,J_2+2,\dots,J_1$,
$
  w_j(\alpha')  \leq C
  2^{-\varepsilon(J_2-J_0)} $.
Using Lemma~\ref{lem:espaleph}, it follows that
\begin{equation}
  \esp \biggl[ \sup_{\alpha'\in(0,2)} |E_2(\alpha')| \biggr]  \leq C
  (J_1-J_2)  2^{-\varepsilon(J_2-J_0)}
   = O(2^{-\xi_2 J})  , \label{eq:bornee2}
\end{equation}
for some $\xi_2>0$ because $\limsup J_0/J < 1/(\alpha+\varepsilon)$. This
concludes the proof.
\end{pf*}

\subsection{Rate of convergence in the stable case}
\label{sec:rate-conv-stable}
\begin{thm} \label{theo:ratecontinustationnaire}
  Let Assumption~\ref{hypo:modele} hold with $\alpha\in(1,2)$ and $p^* = 4$.
  Assume, moreover, that $L_4$ is bounded and that $\mcl(z)= c + O(z^{-\beta})$
  with $c>0$ and $\beta>0$.

  If $\Xgen$ is observed continuously on $[0,T]$, that is, $\bd_{j,k} = d_{j,k}$, then the rate of convergence in probability
  of $\hat\alpha$ is $T^{-\beta /(2\beta+\alpha)}$, obtained for $J_0 = \lfloor
  J/\{2\beta+\alpha\} \rfloor$ and $J_1=J$.

  If $\Xgen$ is observed at discrete time points $1,2,\dots, T$, that is, $\bd_{j,k} = d_{j,k}^D$, then the rate of
  convergence in probability of $\hat \alpha$ is $T^{-\gamma/(2\gamma+\alpha)}$
  with $\gamma = \beta\wedge(2-\alpha)$, obtained for $J_0 = \lfloor
  J/\{2\gamma+\alpha\} \rfloor$ and $J_1=J$.
\end{thm}
\begin{rem}
  Observe that the choice of
  $J_0$ corresponding to the best rate for $\hat\alpha$ depends both on the unknown
  smoothness parameter~$\beta$ and on the parameter $\alpha$ itself. The case of discrete observations is similar to
  that of
  continuous-time observations but with the smoothness parameter $\beta$ replaced by $\gamma= \beta\wedge(2-\alpha)$,
  resulting in a slower rate of convergence.
  This can be explained by the aliasing induced by the interpolation step~(\ref{eq:interpStep}).  It is clear that these
  rates of convergence are the best possible
  for our estimator under the assumption on $\mcl$, since this choice of $J_0$
  makes the squared bias and the variance of the same order of magnitude. However, to our knowledge, the best possible rate of
  convergence for the estimation of $\alpha$ under these observations schemes is an open question. In other words,
  whether our estimator is rate optimal remains unknown.
\end{rem}

The rate of convergence of our estimator is derived under assumptions
on the function $\mathcal L$. The following lemma allows us to check them
through conditions on the joint distribution of $(U,\eta)$.
\begin{lem} \label{lem:rates}
  Let Assumption~\ref{hypo:modele} hold.
\begin{enumerate}[(ii)]
\item[\textup{(i)}\label{item:pareto}]
If there exist positive constants $c$ and $\beta$
  such that, as $t\to\infty$,
\[
L_2(t)=c+O(t^{-\beta}),
\]
then there exists a constant $c'$ such that, as $z\to\infty$,
\begin{equation}\label{eq:rate1}
\mcl (z)= c' +
  \cases{
    O(z^{-\beta}), &\quad if $\beta<2-\alpha$,\cr
    O(z^{\alpha-2}\log z), &\quad if $\beta=2-\alpha$, \cr
    O(z^{\alpha-2}), &\quad if $\beta>2-\alpha$.
}
\end{equation}
\item[(ii)] If there exist positive constants $c$ and $\beta$ such that, as $t\to0$,
\[
\esp[U^2\{1-\cos(\eta t)\}]= c  |t|^{-\alpha} \{ 1  + O(|t|^{\beta}) \}  ,
\]
then there exists a constant $c'$ such that, as $z\to\infty$,
\begin{equation}\label{eq:rate2}
\mcl(z) = c' + O(z^{-\beta})  ,
\end{equation}
provided that $\psi$ belongs to the Sobolev space $W^{(\alpha+\beta)/2-1}$, that is,
\begin{equation}\label{eq:sobolevpsi}
\int_{-\infty}^\infty (1+|\xi|)^{(\alpha+\beta)-2} |\psi^*(\xi)|^2  \d \xi < \infty  ,
\end{equation}
where $\psi^*$ denotes the Fourier transform of $\psi$,
\begin{equation}\label{eq:psiFourier}
\psi^*(\xi)=\int_{0}^M\psi(t)  \rme^{-\rmi\xi t} \d t .
\end{equation}
\end{enumerate}
\end{lem}

\begin{xmpl}
  Assume that $\eta$ has a Pareto distribution, that is, $\pr(\eta>t) = (1 \vee
  t)^{-\alpha}$, and is independent of $U$. This corresponds to
  Lemma~\ref{lem:rates}(i) 
  with $\beta=\infty$, and we can easily
  compute an exact expression for the $O(z^{\alpha-2})$ term:
  \[
    \mcl (z)  =
  c' + \frac{\alpha \esp[U^2]}{2-\alpha} z^{\alpha-2} + o(z^{\alpha-2})  .
  \]
  The best possible rate of convergence of $\hat \alpha$ is thus
  $T^{-(2-\alpha)/(4-\alpha)}$, regardless of the observation scheme.
\end{xmpl}

\begin{xmpl}
  Let $\alpha\in(1,2)$ and suppose that $\eta$ is the absolute value of a
  symmetric $\alpha$-stable random variable. Then Assumption~\ref{hypo:modele}
  holds, say, if $U$ is independent of $\eta$ and has sufficiently many finite
  moments, and
\[
\esp[\cos(\eta t)] = \exp(-\sigma|t|^\alpha)=1-\sigma|t|^\alpha+
O(|t|^{2\alpha}).
\]
By Lemma~\ref{lem:rates}, the best possible rate of convergence of $\hat
\alpha$ is thus $T^{-\gamma/(2\gamma+\alpha)}$ with $\gamma=\alpha$ for
continuous-time observations and $\gamma=2-\alpha$ for discrete-time
observations.
\end{xmpl}

  In the following, we give a decomposition of the error valid under the
  assumption
\[
  0 < \liminf_{T \to \infty} \frac{J_0}{J} \leq \limsup_{T \to \infty} \frac{J_0}{J} < 1.
\]
Optimizing $J_0$ in this decomposition will then give the result.  We use the
same notation as in the proof of Theorem~\ref{theo:stable} with $J_1=J$.  We first give a first rough rate of
convergence for $\hat \alpha$ by adapting the proof of
Theorem~\ref{theo:stable}. Under the present assumptions, $\mcl_0(z)=c$,
which implies $W_0'(\alpha)=0$, and $r^*(z)=O(z^{-\beta})$ as $z\to\infty$.
Then, (\ref{eq:minoration}), (\ref{eq:majo}), (\ref{eq:borneE1})
and (\ref{eq:bornee2}) yield
\begin{equation}\label{eq:rough}
  (\hat \alpha - \alpha)^2 =
  O_P( 2^{-\xi J} + 2^{-\beta J_0}) .
\end{equation}
  Since $\hat\alpha$ is consistent and $\alpha$ is an interior point of the
  parameter set, the first derivative of the contrast function vanishes at
  $\hat\alpha$ with probability tending to one. Hence
  \[
    0 = \frac{\sum_{(j,k)\in\Delta} j  2^{(\hat\alpha-2)j} \bd_{j,k}^2 }
    {\sum_{(j,k)\in\Delta} 2^{(\hat\alpha-2)j} \bd_{j,k}^2 } \ - \ \delta \log(2)  .
  \]
  By the definition of $\delta$, this yields
  \begin{eqnarray}\label{e4.21}
    0 & =& \sum_{(j,k)\in\Delta} (j - \delta)  2^{(\hat\alpha-2)j} \bd_{j,k}^2 \nonumber \\
    & =& \sum_{(j,k)\in\Delta} (j - \delta)  2^{(\alpha-2)j} \bd_{j,k}^2 +
    \log(2) (\hat\alpha - \alpha) \sum_{(j,k)\in\Delta} j (j-\delta)
    2^{(\tilde\alpha-2)j} \bd_{j,k}^2
  \end{eqnarray}
  for a random $\tilde \alpha$ between $\alpha$ and $\hat\alpha$. By the definition
  of $\Lambda_j$, (\ref{e4.21}) implies that
\[
  \hat \alpha - \alpha  = -  \frac{\sum_{j} (j-\delta) 2^{-j} \mathcal
    L(2^j) (1+\Lambda_j)}{\log2\sum_{j} j(j-\delta)
    2^{(\tilde\alpha-\alpha-1)j} \mathcal L(2^j) (1+\Lambda_j)}  .
\]
Denote the sum in the denominator by $D$, and write
\begin{eqnarray*}
  D & = &\sum_{j}j(j-\delta) 2^{-j} \mcl(2^j) + \sum_{j}
  j(j-\delta) \mcl(2^j) 2^{-j}\bigl(2^{(\tilde\alpha-\alpha)j}-1\bigr)(1 + \Lambda_j) \\
  &&{} + \sum_{j} j(j-\delta) 2^{-j} \mcl(2^j) \Lambda_j\\
  & = : & S + R_{1} + R_{2}  .
\end{eqnarray*}
Using Lemma~\ref{lem:svfunc} and~(\ref{eq:delta}), one easily obtains that $S \sim
2^{1-J_0}$ as $J\to\infty$.

Using Lemma~\ref{lem:espaleph}, and the fact that $|\tilde\alpha-\alpha|\leq
|\hat\alpha-\alpha|=o_P(J^{-2})$, one similarly obatins $R_1=o_P(2^{-J_0})$.  To
bound $R_{2}$, we proceed as for bounding $E(\alpha')$ in the proof of
Theorem~\ref{theo:stable} (here with
$\alpha'=\alpha>1$): we write $\sum_j=\sum_{j=J_0+1}^{J_2}+\sum_{j=J_2+1}^{J}$
and apply Lemmas~\ref{lem:espaleph} and~\ref{lem:blob} to obtain
$R_2=o_P(2^{-J_0})$.  Hence, we finally obtain
\begin{equation}
  \hat\alpha-\alpha = \frac{2^{J_0}}{2\log2} \biggl\{ \sum_{j}
    (j-\delta) 2^{-j} \mcl(2^j) \ + \sum_{j} (j-\delta) 2^{-j}
    \mcl(2^j) \Lambda_j \biggr\} \{1+o_p(1)\}  . \label{eq:biaisvariance}
\end{equation}
In~(\ref{eq:biaisvariance}), the terms inside the curly brackets are
interpreted as a deterministic bias term and a stochastic fluctuation term. The
bias is bounded as follows:
\begin{equation}\label{eq:bias}
2^{J_0}\sum_{j} (j-\delta) 2^{-j} \mcl(2^j)  =2^{J_0}\sum_{j} (j-\delta) 2^{-j} (\mcl(2^j)-c)
=O(2^{-\beta J_0})  .
\end{equation}
In the case of continuous-time observations, that is, $\bd_{j,k} = d_{j,k}^S$ or
$\bd_{j,k} = d_{j,k}$, we have
\begin{equation}\label{eq:vartermCase1}
\sum_{j} (j-\delta) 2^{-j} \mcl(2^j) \Lambda_j = O_P\bigl(2^{- J/2 + (\alpha/2-1) J_0 }\bigr) .
\end{equation}
Gathering this bound with~(\ref{eq:biaisvariance}) and~(\ref{eq:bias}), and
setting $J_0 = J/(2\beta+\alpha)$, yields the first claim of
Theorem~\ref{theo:ratecontinustationnaire}, that is,
$\hat\alpha-\alpha=O_P(2^{-\beta/(2\beta+\alpha)})$.

 We now
prove~(\ref{eq:vartermCase1}). Define $\beta_j = n_j^{-1} \sum_{k=0}^{n_j-1}
\{v_j^{-1}(d_{j,k}^{S})^2 -1\}$. Then $\beta_j=\Lambda_j$ if $\bd_{j,k} =
d_{j,k}^S$. Since $\alpha>1$, Lemmas~\ref{lem:varwavcoefstationnaire}
and~\ref{lem:momentbetaj} yield, for some positive constant~$C$,
\begin{equation}\label{eq:betajB}
  \esp[\beta_j^2]=\var(\beta_j)\leq C\frac{L_4(2^j)}{\mcl^2(2^j)} 2^{\alpha
    j-J}  .
\end{equation}
Since $\mcl$ is bounded away from zero and $L_4$ is bounded by assumption, the ratio
$L_4/\mcl^2$ is also bounded.  The Minkowski inequality then yields, for some
constant $C>0$,
\begin{eqnarray}\label{e4.26}
  \esp\biggl[\biggl( \sum_{j} (j-\delta) 2^{-j} \mcl(2^j) \beta_j
\biggr)^2\biggr]^{1/2} &\leq& C 2^{-J/2} \sum_{j} |j-\delta|
  2^{(\alpha/2-1)j}\nonumber\\
  & = & O\bigl(2^{- J/2 + (\alpha/2-1) J_0 }\bigr) .
\end{eqnarray}
If $\bd_{j,k} = d_{j,k}$, we use~(\ref{eq:rjCont}) in Lemma~\ref{lem:espaleph},
and obtain $\esp[|\Lambda_j - \beta_j|] \leq C n_j^{-1/2}$ for some constant
$C>0$. Hence, in this case, since $-1/2 < \alpha/2-1$,
\begin{eqnarray}\label{e4.27}
\esp\biggl[\biggl| \sum_{j} (j-\delta) 2^{-j} \mcl(2^j)
(\Lambda_j -\beta_j) \biggr|\biggr] &\leq& C 2^{-J/2} \sum_{j} |j-\delta|
\mcl(2^j)  2^{-j/2}\nonumber\\
& =& o\bigl(2^{- J/2 + (\alpha/2-1) J_0}\bigr) .
\end{eqnarray}
Inequalities (\ref{e4.26}) and (\ref{e4.27})
imply~(\ref{eq:vartermCase1}).

We now briefly adapt the previous proof to the case of discrete observations.
Define $v_j^D =\esp[(d_{j,k}^{SD})^2]$.  Lemma~\ref{lem:varDiffdjk}(iii)
implies $v_j^D = v_j + O(1)$.  Thus we have $v_j^D =\mcl_D(2^j) 2^{(2-\alpha)j}$ and
\[
\mcl_D(z) = \mcl(z) + O(z^{\alpha-2}) = c + O(z^{-\gamma})  ,
\]
with $\gamma=\beta\wedge(2-\alpha)$. Then, defining
\[
\Lambda_j^D=n_j^{-1} \sum_{k=1}^{n_j} \{(v_j^D)^{-1} \bd_{j,k}^2 -1 \}   ,
\]
we obtain that (\ref{eq:biaisvariance}) still holds with $\mcl_D$ and
$\Lambda^D$ replacing $\mcl$ and $\Lambda$, respectively.
Lemma~\ref{lem:momentbetaDj} implies that $\Lambda_j^D$ has the same order of
magnitude as $\Lambda_j$, so that the stochastic fluctuation term has the same
order of magnitude as in the previous case. The difference comes from the bias
term, which is $O(2^{-\gamma J_0})$. Thus, $\hat\alpha-\alpha = O_P(2^{-\gamma
  J_0}+2^{-\gamma J_0})$, and setting $J_0 = J/(2\gamma+\alpha)$ yields the
second claim of Theorem~\ref{theo:ratecontinustationnaire}.

\section{Concluding remarks}

In this work, we have proved the validity of a wavelet method for the
estimation of the long-memory parameter of an infinite-source Poisson
traffic model, either in a stable or in an unstable state, that is, when
it does or does not converge to a stationary process. We have shown
that a suitable choice of the scales in the estimator (see
Remark~\ref{rem:universalEstimator}) yields a consistent estimator in
both situations, and checked that the estimator is robust to
discrete data sampling.

However, the study of the rates raises some questions concerning the
optimality of this estimator. To draw a comparison, suppose that one
\textit{directly} observes the durations $\eta_1,\dots,\eta_n$ of
clients arriving at times $t_1,\dots,t_n$ in $[0,T]$. Then one can use
the Hill estimator for estimating the tail index $\alpha$. Since $T$
and $n$ are asymptotically proportional and $\eta_1,\dots,\eta_n$ are
independent and identically distributed, the rates of this estimator are those derived in
Hall and Welsh \cite{hallwelsh1984}. In particular, if $\eta$ has a Pareto
distribution, then a parametric rate $\sqrt{T}$ can be obtained. On
the other hand, in the same situation, our wavelet estimator defined
on the observations $\{X(t),t\in[0,T]\}$ has a dramatically
deteriorating rate for $\alpha$ close to 2. It remains to establish
whether this discrepancy comes from the choice of the estimator or
from the fact that the durations $\eta_k$ are not directly observed.

Finally, let us draw a practical conclusion from our study. Care precaution
should be taken with the choice of the scales used in the estimation, as shown by
the conditions on $J_0$ and $J_1$. In particular, if only discrete observations
are available, the best possible rate of convergence is obtained for a much
larger value of $J_0$ than if continous observations are available. Too small a
value of $J_0$ will induce an important bias for finite samples.  Practitioners
should be aware of this restriction and be careful in the interpretation of the
results. These questions will be tackled numerically in a future work.

\begin{appendix}

\section*{Appendix: Technical results}
\label{sec:technical-results}

The following technical lemmas are proved in Fa\"y \textit{et al.} \cite{fayroueffsoulier2005arxiv}.

\begin{lem}  \label{lem:unobs}
  Let Assumption~\ref{hypo:modele} hold.  Let $f$ be a bounded
  measurable compactly supported function such that $\int f(s)  \d s
  = 0$.  Define
\[
  \tilde f (t,v,w) = w \int_t^{t+v} f(s)  \d s   .
\]
Then $\int |\tilde f(t,v,w)|^p   \d t \, \nu(\d v, \d w) < \infty$,
$\esp[N_S( \tilde f)] = 0$ and $\int_0^\infty \Xgen(s) f(s)  \d s
= \break N_S(\tilde f\ind_{\Rset_+\times\Rset_+\times\Rset})$.  If, moreover,
$\esp[\eta]<\infty$, then $N_S(\tilde f) = \int X_S(s)  f(s)  \d s$.
\end{lem}

\begin{lem} \label{lem:momentbetaj}
Let Assumption~\ref{hypo:modele} hold with $p^*\geq4$. Then, there exists a positive constant $C>0$ such that
\[
\var\Biggl(\sum_{k=0}^{n-1} (d_{j,k}^S)^2 \Biggr) \leq C n \bigl\{L_2^2(2^j)
  2^{(4-2\alpha) j} + L_4(2^j)  2^{(3-\alpha) j} \bigr\}  .
\]
Note that the first term dominates for $\alpha<1$ and the second dominates for $\alpha>1$.
\end{lem}

\begin{lem}  \label{lem:unobservableDisc}
  Let Assumption~\ref{hypo:modele} hold. Let $f$ be a bounded measurable
  compactly supported function such that $\int f(s)  \d s = 0$.  Define
\[
  \hat f (t,v,w) = w \int_{-\infty}^\infty g_{t,v}(s) f(s)  \d s
  ,\qquad \check f (t,v,w) = w \int_{-\infty}^\infty h_{t,v}(s) f(s)
  \d s  .
\]
Then, for $p=1,\dots,p^*$, $\int |\hat f(t,v,w)|^p  \d t \, \nu(\d v, \d w) <
\infty$, $\int |\check f(t,v,w)|^p  \d t \, \nu(\d v, \d w) < \infty$, $ \int
\interp[\Xgen] (s)  f(s)  \d s = N_S(\hat
f\ind_{\Rset_+\times\Rset_+\times\Rset})$,
and $\esp[N_S(\hat f)] = \esp[N_S(\check f)] = 0$.  If, moreover,
$\esp[\eta]<\infty$, then $N_S(\hat f) = \int \interp[X_S] (s)  f(s)  \d
s$.
\end{lem}

 Applying Lemma~\ref{lem:unobservableDisc}, we can extend the
 definition of $d_{j,k}^{SD}$ in~(\ref{eq:djk_discrete_case*}) to the
 case $\esp[\eta]=\infty$ by
\begin{equation} \label{eq:defdjksdgeneral}
  d_{j,k}^{SD} = N_S(\hat {\psi}_{j,k}) .
\end{equation}

\begin{lem} \label{lem:varDiffdjk}
\textup{(i)} Let Assumption~\ref{hypo:modele} hold with $p^*\geq1$ and $\alpha\in(0,2)$. Then $\esp
  [d_{j,k}^{SD}] = 0$ for all $j\geq0$ and $k\in\Zset$.

\textup{(ii)} Let
  Assumption~\ref{hypo:modele} hold with $p^*\geq2$ and $\alpha\in(0,2)$.  Then
  $\var(d_{j,k}^S-d_{j,k}^{SD})$ is bounded uniformly for
  $j\in\Nset$ and $k \in\Zset$.

\textup{(iii)}  Let Assumption~\ref{hypo:modele} hold with $p^*\geq2$ and
  $\alpha\in(1,2)$.  Then $|\var(d_{j,k}^S)-\var(d_{j,k}^{SD})|$ is bounded
  uniformly for $j\in\Nset$ and $k \in\Zset$.  \label{item:varsd}
\end{lem}

\begin{lem} \label{lem:espaleph}
  Let Assumption~\ref{hypo:modele} hold with $\alpha \in (1,2)$ and $p^*\geq2$.
Then
\begin{eqnarray}\label{eq:borneespalephj}
  \sup_{0\leq j \leq J} \esp |\Lambda_j| &=& O(1)  ; \\
  \label{eq:rjCont}
\sup_{n\geq1,j\geq0} n^{-1/2}\esp\Biggl[\Biggl|v_j^{-1}\sum_{k=0}^{n-1}\{(d_{j,k}^{S})^2-d_{j,k}^2\} \Biggr|\Biggr]
 & < & \infty  .
\end{eqnarray}
\end{lem}
\begin{lem} \label{lem:blob}
  Let Assumption~\ref{hypo:modele} hold with $\alpha \in (0,2)$ and $p^*\geq4$.
  If $\alpha\leq1/2$, assume (\ref{eq:2ndmomentnul}) and
  (\ref{eq:phiIntLin}).

Let $J^*$ be a sequence depending on $J$ such that $\limsup J^*/J  <
(1/\alpha) \wedge (1/(2-\alpha))$. Then, there exists $\epsilon>0$ such
that
  \begin{equation}
    \sup_{u\in\mcs} \Biggl| \sum_{j=J_0+1}^{J^*} u_j \Lambda_j \Biggr| = O_P(2^{-\epsilon
    J}),
  \end{equation}
where $\mcs$ is the set of sequences $u=(u_0,\dots)$ satisfying $\sum_{j \in \Nset} |u_j| \leq 1$.
\end{lem}

\begin{lem} \label{lem:momentbetaDj}
  Let Assumption~\ref{hypo:modele} hold with $p^*\geq4$ and
$\alpha\in(1,2)$. Then, there exists a positive constant $C>0$ such
that
  \begin{equation} \label{eq:varsumsquarestationnaire}
    \var\Biggl(\sum_{k=0}^{n-1} (d_{j,k}^{SD})^2 \Biggr) \leq C L_4(2^j) n
    2^{(3-\alpha)j} .
  \end{equation}
\end{lem}

\begin{lem} \label{lem:svfunc}
Let $\rho$ be a positive real and $\rho':= (2^{\rho} -
  1)^{-1}$.   Let $\ell^*$ be a non-increasing function on $[1,\infty)$ such that
  $\lim_{s\to\infty} \ell^*(s) = 0$, and let $\ell$ be a function on
  $[1,\infty)$ such that $|\ell(s)| \leq \ell^*(s)$ for all
  $s\in[1,\infty)$.  Define
  \[
    L(x) = c \exp \biggl\{ \int_{1}^x \frac{\ell(s)}s  \d
    s\biggr\}\quad
        \mbox{and}\quad
    \omega_j = \frac{2^{-\rho j} L(2^j)}{\sum_{j'=J_0+1}^{J_1} 2^{-\rho j'}
      L(2^{j'})}  .
  \]
Then, as $J_0 \to \infty$ and for any $\epsilon > 0$,
  \begin{eqnarray}
    \sum_{j=J_0+1}^{J_1} \omega_j j  &=& J_0 + 1 + \rho'\bigl(1+ O ( \ell^*(2^{J_0}))\bigr) +
    O\bigl(J_1(2-\epsilon)^{J_0 - J_1}\bigr), \label{eq:firstorder}\\
  \label{eq:secondorder}
    \sum_{j=J_0+1}^{J_1} \omega_j j^2 &=&  J_0^2 + 2J_0(1+\rho')
  + 2\rho'^2 +3\rho' + 1   + \rho' O ( \ell^*(2^{J_0}))\nonumber
\\
&&{}  +
    O\bigl(J_1^2(2-\epsilon)^{J_0 - J_1}\bigr).
  \end{eqnarray}
\end{lem}
\end{appendix}

\printhistory

\end{document}